\documentclass[12pt]{amsart}
\usepackage{amscd,amsmath,amsthm,amssymb}
\usepackage{pstcol,pst-plot,pst-3d}
\usepackage{lmodern,pst-node}
\def\frk{\frak}               

\def\Phi{{\frk n}}
\def\Phi{{\frk N}}
\def\MI{{\mathcal I}}
\def\MJ{{\mathcal J}}
\def\ML{{\mathcal L}}
\def\MR{{\mathcal R}}
\def\Ic{{\mathcal I}}

\def\Lc{{\mathcal L}}
\def\opn#1#2{\def#1{\operatorname{#2}}} 
\opn\chara{char} \opn\length{\ell} \opn\pd{pd} \opn\rk{rk}
\opn\projdim{proj\,dim} \opn\injdim{inj\,dim} \opn\rank{rank}
\opn\depth{depth} \opn\grade{grade} \opn\height{height}
\opn\embdim{emb\,dim} \opn\codim{codim}

\opn\Tr{Tr} \opn\bigrank{big\,rank}
\opn\superheight{superheight}\opn\lcm{lcm}
\opn\trdeg{tr\,deg}
\opn\reg{reg} \opn\lreg{lreg} \opn\ini{in} \opn\lpd{lpd}
\opn\size{size} \opn\sdepth{sdepth}
\opn\link{link}\opn\fdepth{fdepth}\opn\lex{lex}
\opn\div{div} \opn\Div{Div} \opn\cl{cl} \opn\Cl{Cl}
\opn\Spec{Spec} \opn\Supp{Supp} \opn\supp{supp} \opn\Sing{Sing}
\opn\Ass{Ass} \opn\Min{Min}\opn\Mon{Mon}
\opn\Ann{Ann} \opn\Rad{Rad} \opn\Soc{Soc}
\opn\Im{Im} \opn\Ker{Ker} \opn\Coker{Coker} \opn\Am{Am}
\opn\Hom{Hom} \opn\Tor{Tor} \opn\Ext{Ext} \opn\End{End}
\opn\Aut{Aut} \opn\id{id}

\opn\nat{nat}
\opn\pff{pf}
\opn\Pf{Pf} \opn\GL{GL} \opn\SL{SL} \opn\mod{mod} \opn\ord{ord}
\opn\Gin{Gin} \opn\Hilb{Hilb}\opn\sort{sort}
\opn\aff{aff} \opn\con{conv} \opn\relint{relint} \opn\st{st}
\opn\lk{lk} \opn\cn{cn} \opn\core{core} \opn\vol{vol}
\opn\link{link} \opn\star{star}\opn\lex{lex}\opn\set{set}
\opn\gr{gr}
\def\Rees{{\mathcal R}}
\def\pot#1#2{#1[\kern-0.28ex[#2]\kern-0.28ex]}

\opn\dirlim{\underrightarrow{\lim}}
\opn\inivlim{\underleftarrow{\lim}}
\let\union=\cup
\let\sect=\cap

\let\iso=\cong
\let\Union=\bigcup
\let\Sect=\bigcap

\let\to=\rightarrow

\def\Implies{\ifmmode\Longrightarrow \else
        \unskip${}\Longrightarrow{}$\ignorespaces\fi}
\def\implies{\ifmmode\Rightarrow \else
        \unskip${}\Rightarrow{}$\ignorespaces\fi}
\def\iff{\ifmmode\Longleftrightarrow \else
        \unskip${}\Longleftrightarrow{}$\ignorespaces\fi}

\let\:=\colon
\newtheorem{Theorem}{Theorem}[section]
\newtheorem{Lemma}[Theorem]{Lemma}
\newtheorem{Corollary}[Theorem]{Corollary}
\newtheorem{Proposition}[Theorem]{Proposition}

\let\epsilon\varepsilon
\let\kappa=\varkappa
\def\qed{\ifhmode\textqed\fi
      \ifmmode\ifinner\quad\qedsymbol\else\dispqed\fi\fi}
\def\textqed{\unskip\nobreak\penalty50
       \hskip2em\hbox{}\nobreak\hfil\qedsymbol
       \parfillskip=0pt \finalhyphendemerits=0}
\def\dispqed{\rlap{\qquad\qedsymbol}}

\opn\dis{dis}
\def\pnt{{\raise0.5mm\hbox{\large\bf.}}}

\opn\Lex{Lex}

\begin{document}

\title {Monomial ideals and toric rings of Hibi type arising  from  a finite poset}

\author {Viviana Ene, J\"urgen Herzog and Fatemeh Mohammadi}

\address{Faculty of Mathematics and Computer Science, Ovidius University, Bd.\ Mamaia 124,
 900527 Constanta, Romania} \email{vivian@univ-ovidius.ro}

\address{J\"urgen Herzog, Fachbereich Mathematik, Universit\"at Duisburg-Essen, Campus Essen, 45117
Essen, Germany} \email{juergen.herzog@uni-essen.de}

\address{Fatemeh Mohammadi, Department of Pure Mathematics, Ferdowsi University of Mashhad
P.O. Box 1159-91775, Mashhad, Iran}\email{fatemeh.mohammadi716@gmail.com}

\begin{abstract}
In this paper we study monomial ideals attached to posets,  introduce generalized Hibi rings and investigate their algebraic and homological properties. The main tools to study these objects are Gr\"obner basis theory, the concept of sortability due to Sturmfels and the theory of weakly polymatroidal ideals.
\end{abstract}

\maketitle

\section*{Introduction}

In 1985  Hibi \cite{Hibi1} introduced a class of algebras which nowadays  are called Hibi rings. They are  toric rings attached to finite posets, and may be viewed as natural generalizations of polynomial rings. Indeed, a polynomial ring  in $n$ variables  over a field $K$ is just the  Hibi ring of the poset $[n]=\{1,2,\ldots,n\}$.

Hibi rings appear naturally in various combinatorial and algebraic contexts, for example in invariant theory. Hodge algebras may be viewed as flat deformations of Hibi rings. In this sense the  coordinate ring of the flag variety for $\GL_n$  is the deformation of the Hibi ring for the so-called  Gelfand-Tsetlin poset.

Given a finite poset $P=\{p_1,\ldots,p_n\}$, let  $\MI(P)$ be the ideal lattice of $P$. By Birkhoff's theorem any finite distributive lattice arises in this way. Let $K$ be field. Then the Hibi ring  over $K$ attached to  $P$ is the toric ring $K[\MI(P)]$ generated by the set  of  monomials $\{x_It\:\; I\in \MI(P)\}$ where $x_I=\prod_{p_i\in I}x_i$. Let $T=K[\{y_I\:\; y_I\in \MI(P)\}]$ be the polynomial ring in the variables $y_I$ over $K$, and $\varphi\:T\to K[\MI(P)]$ the $K$-algebra homomorphism with $y_I\mapsto x_It$. One fundamental result concerning Hibi rings is that the toric ideal $L_P=\Ker \varphi$ has a reduced Gr\"obner basis consisting of the so-called Hibi relations:
\[
y_Iy_I-y_{I\sect J}y_{I\union J}\quad \text{with} \quad I\not\subseteq J\quad \text{and}\quad J\not\subseteq I.
\]
Hibi showed \cite[Theorem 2.c]{Hibi1} that any Hibi ring is a normal Cohen--Macaulay domain, and that it is Gorenstein if and only if the attached poset $P$ is pure  \cite[Corollary 3.d]{Hibi1}, that is, all maximal chains of $P$  have the same cardinality.

More generally,  for any lattice $\ML$, not necessarily distributive, one may consider the $K$ algebra $K[\ML]$ with generators $y_\alpha$, $\alpha\in \ML$, and relations $y_\alpha y_\beta=y_{\alpha\wedge \beta}y_{\alpha\vee \beta}$ where $\wedge$ and $\vee$ denote meet and join in $\ML$.  Hibi showed that $K[\ML]$ is a domain if and only if $\ML$ is distributive, in other words,  if $\ML$ is an ideal lattice of a poset.

The starting point of this paper are the so-called Hibi ideals  which were first introduced in \cite{cm}. Attached to each finite poset $P=\{p_1,\ldots,p_n\}$, one defines the Hibi ideal $H_P$ as the monomial ideal in the polynomial ring  $K[x_1,\ldots,x_n,y_1,\ldots,y_n]$ generated by the monomials  $x_Iy_{P\setminus I}$ with $I\in \MI(P)$. Note that the toric ring generated over $K$ by these monomials is isomorphic to the Hibi ring attached to $P$. The significance of Hibi ideals is that their Alexander dual can be interpreted as the edge ideal of a bipartite graph. To be precise, if we define the bipartite graph $G$ on the vertex set $\{x_1,\ldots,x_n,y_1,\ldots,y_n\}$  by saying that $\{x_i,y_j\}$ is an edge of $G$ if and only $p_i\leq p_j$, then $H_P^\vee$ is the edge ideal of $G$ in the sense of Villarreal \cite{Vi}. It turned out  that bipartite graphs obtained in this way are exactly the Cohen--Macaulay  bipartite graphs.

Motivated by the dual relationship between Hibi ideals and edge ideals of bipartite graphs we introduce in this paper the following ideals attached to a finite poset $P=\{p_1,\ldots,p_n\}$: fix  integers $r\geq s\geq 1$, a  field  $K$  and let $S$ be the polynomial ring over $K$ in the variables $x_{ij}$ with  $i=1,\ldots, r$ and $j=1,\ldots,n$.

We denote by $\mathcal{C}_r(P)$ the set of   multichains  of length $r$. In other words, the elements of $\mathcal{C}_r(P)$ are subsets
$\{p_{j_1},\ldots, p_{j_r}\}$ of $P$ with $p_{j_1}\leq p_{j_2}\leq \cdots \leq p_{j_r}$.

For $C\in \mathcal{C}_r(P)$,  $C=\{p_{j_1},\ldots, p_{j_r}\}$ and $\emptyset \neq S\subseteq [r]$ we set
\[
u_{C,S} =\prod_{i\in S}x_{ij_i},
\]
and  $u_C=u_{C,[r]}$. Then we define the monomial ideals
\[
I_{r,s}(P)=(u_{C,S}\:\; C\in \mathcal{C}_r(P)\quad \text{and}\quad   S\subset [r] \quad \text{with}\quad |S|=s),
\]
and $H_{r,s}(P)=I_{r,s}(P)^\vee$.

We  call $I_{r,s}(P)$ the {\em multichain ideal} of type $(r,s)$,  and $H_{r,s}(P)$  the {\em generalized Hibi ideal}  of type $(r,s)$ of the poset $P$, since  $H_{2,2}(P)$ is just the classical Hibi ideal $H_P$. For simplicity the  ideals $H_{r,r}(P)$ will  be denoted by $H_r(P)$. It is worthwhile to notice that the ideals $I_{r,s}(P)$ may be interpreted as facet ideals of completely balanced simplicial complexes, as introduced by Stanley \cite{Sta}.

In Theorem~\ref{FatemehViviana} we compute explicitly the minimal monomial set of generators of $H_{r}(P)$ and show that $H_{r,s}(P)=H_r(P)^{\langle r-s+1\rangle}$, where $I^{\langle k\rangle}$ denotes the $k$th squarefree power of a  squarefree monomial ideal. It turns out that the generators of $H_r(P)$ correspond bijectively to chains of length $r$
\[
\MI\: I_1\subseteq I_2\subseteq \cdots \subseteq I_r=P
\]
of poset ideals of $P$.

Based on this explicit description we show in Theorem~\ref{power} that all powers $H_r(P)^k$ of the ideal $H_r(P)$ are weakly polymatroidal. The concept of weakly polymatroidal ideals has been introduced by Hibi and Kokubo in \cite{Hibi2}
where they showed that they share  with polymatroidal ideals  the nice property of having linear quotients. In particular, we conclude from this  that  the ideals  $H_r(P)^k$ have a linear resolution for all $k\geq 1$. Now we   use the fact shown in Theorem~\ref{FatemehViviana} that $H_{r,s}(P)$ is a suitable squarefree power of the ideal $H_r(P)$, and observe that the  squarefree  part of a weakly polymatroidal ideal is again weakly polymatroidal (see Lemma~\ref{squarefree}) to  finally deduce in Theorem~\ref{irscm} that all Hibi ideals $H_{r,s}(P)$ are weakly polymatroidal. By the Eagon--Reiner Theorem \cite{ER} this implies that  for any  finite poset and all integers $1\leq s\leq r$ the chain ideals $I_{r,s}(P)$ are  Cohen--Macaulay. Thus we obtain a rich family of completely balanced simplical complexes whose facet ideals are  Cohen--Macaulay. In the case that $s=2$ this yields  a class of $r$-partite graphs with Cohen--Macaulay edge ideal.

In Corollary~\ref{ircm} we show that the ideals $I_r(P)$ are Gorenstein,  if and only if they are generated by a regular sequence which is the case if and only all elements of $P$ are pairwise incomparable.

Section 3 is devoted to the study of the resolution of the ideal $H_r(P)$. As $H_r(P)$ has linear quotients, the resolution can in principle be obtained as an iterated mapping cone. To get an explicit description of the maps in the resolution one has to know all the linear quotients. This is described in Lemma~\ref{sets}. With this  information at hand we can describe the projective dimension of $H_r(P)$ as the
maximal length of antichains in $P$, see Corollary~\ref{H_r(P)}. Applying then a result of Terai which relates the projective dimension of an ideal to the regularity of its Alexander dual we obtain a nice formula for the regularity for  the chain ideals $I_r(P)$. Next we show that the ideals $H_r(P)$ have regular decomposition functions in  the sense of \cite{HT}, and then apply  a result of the same paper to finally obtain in Theorem~\ref{resolution} the explicit resolution of the ideals $H_r(P)$.

The remaining sections of the paper are devoted to the study of the toric rings $R_{r,s}(P)$ which naturally generalize the classical  Hibi rings. The toric ring $R_{r,s}(P)$, respectively $R_{r}(P)$,  is defined to be the standard graded $K$-algebra generated over $K$ by the unique minimal minimal  set $G_{r,s}(P)$ of monomial generators of $H_{r,s}(P)$, respectively of $H_{r}(P)$. In Theorem~\ref{generalized} we first  show that  $R_{r}(P)$ has a quadratic reduced Gr\"obner basis consisting of Hibi type relations. This result is used to show in Corollary~\ref{normalcm} that  all the toric rings $R_r(P)$ are normal Cohen--Macaulay domains and to identify in Theorem~\ref{hibiring} the toric ring $R_r(P)$ as a classical Hibi ring attached to the direct product $P\times Q_{r-1}$ of the  poset $P$ and the poset $Q_{r-1}=[r-1]$.  By using  this fact  we see in  Corollary~\ref{hibiringgor} that  $R_r(P)$ is  Gorenstein if and only if $P$ is pure.

The situation for the toric rings $R_{r,s}(P)$ is more complicated. Among their relations are also relations which are not of Hibi type and so these algebras cannot be identified with classical Hibi rings for  suitable posets. Indeed, if we choose for $P$ the poset consisting just of one element, then $H_{r,s}(P)$ is nothing than the squarefree Veronese $K$-algebra which is generated over $K$ by all squarefree monomials of degree $s$ in $r$ variables. If we choose the same monomial order to compute the Gr\"obner basis of the corresponding toric ideal of this algebra as we did in the proof for the algebras $R_r(P)$, then in this particular  case this term order is just the reverse lexicographic order induced by the order of the variables which is given by the lexicographic order of the generators of the algebra.  For this monomial order the algebra $R_{r,s}(P)$  has  a reduced Gr\"obner basis consisting of binomials of degree $\leq 3$ with squarefree initial ideal (Theorem~\ref{veroneseviv}).   That some of the elements in the reduced Gr\"obner basis may  indeed be of degree $3$,  can be seen for example if we choose $r=6$ and $s=3$.

The question arises whether there is a monomial order for which the algebras $R_{r,s}(P)$ has a quadratic Gr\"obner basis. The answer is yes, and the method to prove this is due to Sturmfels who used a sorting order  to show that all algebras of Veronese type  have a quadratic Gr\"obner basis. What we need to show is that  the set of monomial $G_{r,s}(P)$ is sortable. This then implies by a theorem of Sturmfels \cite[Theorem 14.2]{St} that the quadratic sorting relations with the unsorted pairs as leading terms form a  Gr\"obner basis with respect to the sorting order induced by the sorting of the monomials.   We prove in Theorem~\ref{fatemehsortable} that the set $G_{r,s}(P)$ is indeed sortable. As a consequence  we obtain  that the algebras $R_{r,s}(P)$ are again all normal Cohen--Macaulay domains.

In the last section we study the Rees algebra of the ideals $H_{r,s}(P)$. In \cite{HHM} the authors introduce the so-called $\ell$-exchange property  which guarantees that the Rees algebra of a monomial ideal which is generated in one degree has a Gr\"obner basis composed of the Gr\"obner basis of the fibre of the Rees algebra and binomial relations which are  linear in the variables of the base ring. It is shown in Proposition~\ref{exchange} that the $\ell$-exchange property is satisfied  for sortable, weakly polymatroidal ideals. Thus we may apply the result of \cite{HHM} and find that the Rees algebra $\mathcal{R}(H_{r,s}(P))$ has a quadratic Gr\"obner basis. As applications we find that all powers of $H_{r,s}(P)$ are normal and have  a linear resolution,  and that $\mathcal{R}(H_{r,s}(P))$  is a normal Cohen--Macaulay Koszul algebra.

While we can  characterize the Gorenstein ideals   $I_r(P)$ and the Gorenstein rings $R_r(P)$, we do not have such a characterization for the ideals
$I_{r,s}(P)$ and the rings $R_{r,s}(P)$. 

\section{Multichain ideals of a poset  and their Alexander dual}

In this section we determine the Alexander dual $H_{r,s}(P)$ of the multichain ideal $I_{r,s}(P)$. Recall from the introduction that
\[
I_{r,s}(P)=(u_{C,S}\:\; C\in \mathcal{C}_r(P)\quad \text{and}\quad   S\subset [r] \quad \text{with}\quad |S|=s),
\]
where $\mathcal{C}_r(P)$ is the set of multichains of length $r$ in $P$,  and where
\[
u_{C,S} =\prod_{i\in S}x_{ij_i}.
\]

In order to formulate the main result of this section we introduce some notation. Given a multichain
\[
\mathcal{I}\:\; I_1\subseteq I_2\subseteq \ldots \subseteq I_r=P
\]
of poset ideals in $P$, we attach to it the following monomial in $S$:
\[
u_\MI=x_{1J_1}x_{2J_2}\cdots x_{rJ_r},\quad \text{where $x_{jJ_j}=\prod_{p_k\in J_j}x_{jk}$ and $J_j=I_j\setminus I_{j-1}$\quad for  $j=1,\ldots,r$.}
\]
We denote by $\mathcal{I}_r(P)$ the set of multichains of poset ideals of length $r$ in $P$, and for any squarefree monomial ideal $L$ we denote by $L^{\langle k\rangle}$ the $k$th squarefree power of $L$, that is, the ideal generated by all squarefree elements in $L^k$.

\begin{Theorem}
\label{FatemehViviana}
Let $P$ be a finite poset. Then
\begin{enumerate}
\item[(a)] The Alexander dual $H_r(P)$ of $I_r(P)$ is the ideal $(u_{\MI}\:\; \MI\in \mathcal{I}_r(P))$;
\item[(b)] The Alexander dual  $H_{r,s}(P)$ of $I_{r,s}(P)$ is the ideal $H_r(P)^{\langle r-s+1\rangle}$.
\end{enumerate}
\end{Theorem}

\begin{proof}
First we show
that for any multichain of poset elements
$p_{\ell_1}\leq \cdots \leq p_{\ell_r}$, the ideal \[Q=(x_{1\ell_1},\ldots,x_{r\ell_r})\] is a minimal prime ideal of  $H_{r}(P)$.

In order to see that $Q$ contains $H_{r}(P)$, we show that for each  $u_{\MI}\in H_{r}(P)$  there exists some $j\in [r]$ such that $x_{j,\ell_{j}}$ divides $u_{\MI}$.

By contrary assume that no $x_{j,\ell_{j}}$ divides $u_{\MI}$. Then
$p_{\ell_r}\not\in I_r\setminus I_{r-1}$, and so $p_{\ell_r}\in I_{r-1}$. Since $p_{\ell_{r-1}}\leq p_{\ell_r}$ it follows that $p_{\ell_{r-1}}\in I_{r-1}$. On the other hand, since  $p_{\ell_{r-1}}\not\in I_{r-1}\setminus I_{r-2}$ we have $p_{\ell_{r-1}}\in I_{r-2}$ which implies that  $p_{\ell_{r-2}}\in I_{r-2}$. Continuing in this way we get $p_{\ell_1}\in I_1$, a contradiction.

Suppose $Q$ is not a minimal prime ideal, then there exists an integer $i$ such that $Q_i=(x_{1\ell_1},\ldots,x_{i-1\ell_{i-1}},x_{i+1\ell_{i+1}}, \ldots, x_{r\ell_r})$ contains $H_{r}(P)$. However  $x_{i P}$ is a generator of $I_r(P)$  which does not belong to $Q_i$, a contradiction.

Next we show that for any multichain of poset elements
$p_{\ell_1}\leq \cdots \leq p_{\ell_r}$ and each subset  $\{\ell_{t_{1}}, \ldots, \ell_{t_{r-s}}\}$ of $\{\ell_{1}, \ldots, \ell_r\}$, the ideal

\begin{eqnarray}
\label{primeideal}
Q=(\{x_{1\ell_1},\ldots,x_{r\ell_r}\}\setminus \{x_{{t_1}\ell_{t_1}},\ldots,x_{{t_{r-s}}\ell_{t_{r-s}}}\})
\end{eqnarray}
is a prime ideal containing $H_{r}(P)^{\langle r-s+1\rangle}$.

Let $u=u_{\MI_1}\cdots u_{\MI_{r-s+1}}$ be an arbitrary element in $H_{r}(P)^{\langle r-s+1\rangle}$.
We show that there exists some $j\in [r]\setminus \{t_1,\ldots,t_{r-s}\}$ such that $x_{j,\ell_{j}}$ divides $u$.

We know already that $(x_{1\ell_1},\ldots,x_{r\ell_r})$ is a minimal prime ideal of $H_r(P)$. So for each $1\leq k\leq r-s+1$, there exists an index  $j_k$ such that $x_{j_k,\ell_{j_k}}|u_{{\MI}_k}$. Since $u$ is a squarefree monomial, the elements  $x_{j_1,\ell_{j_1}},\ldots,x_{j_{r-s+1},\ell_{j_{r-s+1}}}$  are pairwise distinct, and hence at least one of them  must belong to $\{x_{1\ell_1},\ldots,x_{r\ell_r}\}\setminus \{x_{{t_1}\ell_{t_1}},\ldots,x_{{t_{r-s}}\ell_{t_{r-s}}}\}$, as we wanted to show.

In order to show that $Q$ is a minimal prime ideal of $H_{r}(P)^{\langle r-s+1\rangle}$ we first observe the following fact: for each minimal prime ideal $Q'$ of $H_{r}(P)^{\langle r-s+1\rangle}$, there exist $s$ indices $j_1<j_2<\cdots < j_{s}$ in $[r]$ such that
\[
Q'\cap \{x_{j_i1},\ldots,x_{j_in}\}\neq \emptyset \quad \text{for} \quad i=1,\ldots,s.
\]
By contrary, there exist $t_1<t_2<\cdots <t_{r-s+1}$ with $t_i\in[r]$ such that
\[
Q'\cap \{x_{t_i1},\ldots,x_{t_in}\}=\emptyset \quad \text{for}\quad i=1,\ldots,r-s+1.
\]
Then the monomials $u=\prod_{i=1}^{r-s+1}x_{t_iP}\in H_{r}(P)^{\langle r-s+1\rangle}$, but $u\not\in Q'$, a contradiction.

It follows from these considerations that each minimal prime ideal of $H_{r}(P)^{\langle r-s+1\rangle}$ has at least $s$ variables as generators. Since $Q$ has precisely $s$ variables as generators, the prime ideal $Q$ must be a minimal prime ideal of $H_{r}(P)^{\langle r-s+1\rangle}$.

\medskip
It remains to be shown that each minimal prime of $H_{r}(P)^{\langle r-s+1\rangle}$ is of the form  (\ref{primeideal}). So let $Q$ be an arbitrary minimal prime ideal of  $H_{r}(P)^{\langle r-s+1\rangle}$. We know from the proof before  that $Q$ has to have exactly $s$ variables as generators. Assume that $Q=(x_{i_1\ell_{1}}, x_{i_2\ell_{2}}, \ldots, x_{i_{s}\ell_{{s}}})$ for some $i_1<i_2<\cdots<i_{s}$. We are going to show that $p_{l_1}\leq p_{l_2}\leq\cdots\leq p_{l_{s}}$. By contrary, assume that for some $j$, $p_{\ell_j}\nleq p_{\ell_{j+1}}$. Then consider the multichain $\MI:\ I_{1}\subseteq I_2\subseteq \cdots\subseteq I_r$ of poset ideals of $P$ with
\[
I_{1}=\cdots=I_{i_{j-1}}=\emptyset,\quad I_{i_j}=\cdots=I_{i_{j+1}-1}=\langle p_{\ell_{j+1}}\rangle,\quad \text{and}\quad I_{i_{j+1}}=\cdots=I_r=P,
\]
where for $p\in P$ we set $\langle p\rangle =\{q\leq p\:\; q\in P\}$. Therefore, $u_{\MI}=x_{i_j\langle p_{\ell_{j+1}}\rangle}x_{i_{j+1}P\setminus \langle p_{\ell_{j+1}}\rangle}$.

\medskip
Let $\{t_1,t_2,\ldots,t_{r-s}\}=[r]\setminus\{i_1,i_2,\ldots,i_{s}\}$, and let
$u=u_{\MI} \prod_{i=1}^{r-s}x_{t_iP}$. Then
$u\in H_{r}(P)^{\langle r-s+1\rangle}$, but $u\not\in Q$, a contradiction.
\end{proof}

\section{Generalized Hibi ideals and their powers}

Kokubo  and Hibi in \cite{Hibi2} introduced weakly polymatroidal
ideals as a generalization of polymatroidal ideals. They show \cite[Theorem~1.4]{Hibi2} that weakly polymatroidal ideals have linear quotients. In particular, this implies that weakly polymatroidal ideals have a linear resolution.

Let $R=K[x_1,\ldots,x_n]$ be a polynomial ring over the field $K$.
Recall that a monomial ideal $I\subset R$ which is generated in one degree
is called {\em weakly polymatroidal} if for any two monomials
$u=x_1^{a_1}\cdots x_n^{a_n}$ and $v=x_{1}^{b_1}\cdots x_n^{b_n}$ in
$G(I)$ for which there exists an integer $t$ with   $a_1=b_1,\ldots,a_{t-1}=b_{t-1}$ and $a_t>b_t$,
there exists $\ell>t$ such that $x_t(v/x_{\ell})\in I$. Here we denote as usual the unique minimal set of monomial generators of $I$. Note that the concept  weakly polymatroidal depends on the order of the variables of $R$.

Observe that a squarefree monomial ideal $I$ which is generated in one degree is weakly polymatroidal if for any  two monomials
$u=x_{i_1}\cdots x_{i_d}$ with $i_1<i_2<\cdots <i_d$, and $v=x_{j_1}\cdots x_{j_d}$ with $j_1<j_2<\cdots <j_d$ in
$G(I)$ such that $i_1=j_1,\ldots,i_{t-1}=j_{t-1}$ and $i_t<j_t$,
there exists $\ell\geq t$ such that $x_{i_t}(v/x_{j_{\ell}})\in I$.

\medskip
We define a partial order on the set $\mathcal{I}_r(P)$ by setting $\MI\leq \MI'$ if $I_i\subseteq  I_i'$ for $i=1,\ldots,r$. Observe that the partially ordered set $\mathcal{I}_r(P)$ is a distributive lattice, if we define the meet of $\MI\: I_1\subseteq \cdots \subseteq I_r$  and $\MI'\: I'_1\subseteq \cdots \subseteq I'_r$  as $\MI\sect \MI'$ where  $(\MI\sect \MI')_i=I_i\sect I_i'$ for $i=1,\ldots,r$, and the join as $\MI\union \MI'$ where $(\MI\union \MI')_i=I_i\union I'_i$ for $i=1,\ldots,r$,

The following lemma was
proved in \cite[p.99]{Hibi1} for $r=2$ and $k=1$.

\begin{Lemma}\label{describe}
Any element in the minimal generating set of $H_r(P)^{k}$ can be written as ${u}_ {\MI_1}
\cdots u_ {\MI_k}$, where $\MI_i\in \mathcal{I}_r(P)$ and ${\MI_1} \leq\cdots
\leq \MI_k$.
\end{Lemma}

\begin{proof}
We first  claim that
\begin{eqnarray}
\label{unionsect}
u_{\MI}u_{{\MI}'}={u}_ {{\MI}\cap {\MI}'}{u}_ {{\MI}\cup {\MI}'}.
\end{eqnarray}
Indeed, the equation (\ref{unionsect}) is valid if and only if
\[
\frac{x_{t I_t}}{x_{t I_{t-1}}}\cdot \frac{x_{t I'_t}}{x_{t I'_{t-1}}}= \frac{x_{t I_t\sect I_t'}}{x_{t I_{t-1}\sect I'_{t-1}}}\cdot
\frac{x_{t I_t\union I_t'}}{x_{t I_{t-1}\union I'_{t-1}}}
\]
for $t=1,\ldots r$.

In order to see that this identity holds, just observe  that
\[
x_{t I_t\sect I_t'}= \gcd\{x_{t I_t}, x_{t I_t'}\},\quad   x_{t I_{t-1}\sect I_{t-1}'}= \gcd\{x_{t I_{t-1}}, x_{t I_{t-1}'}\},
\]
and
\[
x_{t I_t\union  I_t'}=\frac{x_{t I_t}\cdot x_{t I_t'}} { \gcd\{x_{t I_t}, x_{t I_t'}\}}\quad x_{t I_{t-1}\union  I_{t-1}'}=
\frac{x_{t I_{t-1}}\cdot x_{t I_{t-1}'}} { \gcd\{x_{t I_{t-1}}, x_{t I_{t-1}'}\}}.
\]
\medskip

Now let $u=u_{\MJ_1}\cdots u_{\MJ_k}\in H_r(P)^{\langle k\rangle}$. By induction we show that $u$ can be written as
$u_{\MJ'_1}\cdots u_{\MJ'_k}$ such that ${\MJ'_i}\leq \MJ'_k$ for
$i=1,\ldots,k-1$. Indeed, applying the induction hypothesis we may assume  $u_{\MJ_1}\cdots u_{\MJ_{k-1}}=u_{\MJ'_1}\cdots u_{\MJ'_{k-1}}$ with ${\MJ'_i} \leq \MJ'_{k-1}$ for $i=1,\ldots,k-2$. Then, by using (\ref{unionsect}) we obtain
\begin{eqnarray*}
u_{\MJ_1}\cdots u_{\MJ_k}
&=&(u_{\MJ'_1}\cdots {u}_{\MJ'_{k-2}})(u_{\MJ'_{k-1}} {u}_{\MJ_{k}})
\\
&=&({u}_{\MJ'_1}\cdots {u}_{\MJ'_{k-2}})(u_{\MJ'_{k-1}\cap \MJ_k}u_{\MJ'_{k-1}\cup{\MJ_{k}}}).
\end{eqnarray*}
We have $\MJ'_{k-1}\cap \MJ_k\leq {\MJ'_{k-1}\cup{\MJ_{k}}}$ and $\MJ'_i\leq
{\MJ'_{k-1}\cup{\MJ_{k}}}$ for $i=1,\ldots, k-2$.
Hence we may indeed assume from the very beginning that in $u=u_{\MJ_1}\cdots u_{\MJ_k}$, we have
${\MJ_i} \leq \MJ_k$ for
 $i=1,\ldots,k-1$.

Now we show that one can write $u_{\MJ_1}\cdots u_{\MJ_k}={u}_ {\MI_1} \cdots
{u}_ {\MI_{k}}$ with ${\MI_1} \leq \cdots \leq \MI_k$.
By induction assume that ${u}_ {\MJ_1} \cdots
{u}_ {\MJ_{k-1}}={u}_ {\MI_1} \cdots
{u}_ {\MI_{k-1}}$ with ${\MI_1}
\leq \cdots \leq \MI_{k-1}$. Since
$\bigcup_{i=1}^{k-1}\MI_i=\bigcup_{i=1}^{k-1}\MJ_i\leq \MJ_k$, setting
$\MI_k=\MJ_k$ we have ${\MI_1} \leq \cdots \leq \MI_k$ and
$u_{\MJ_1}\cdots u_{\MJ_k}=u_{\MI_1}\cdots u_{\MI_k}$.
\end{proof}

\begin{Theorem}\label{power}
For any positive integer $k$, the ideal $H_r(P)^{k}$ is weakly polymatroidal.
\end{Theorem}

\begin{proof}
Let $P=\{p_1,\ldots,p_n\}$. We may assume that if  $p_i< p_j$, then  $i<j$. We are going to show that  $H_r(P)^{k}$ is weakly polymatroidal with respect to the following order
\[
x_{11}>x_{12}>\cdots>x_{1n}>x_{21}>\cdots>x_{2n}>\cdots>x_{r1}>\cdots>x_{rn}
\]
of the variables.

\medskip

Let $u=u_{\MI_1}\cdots
u_{\MI_k}$ with ${\MI_k}
\leq \cdots \leq \MI_1$ and $v=u_{\MJ_1}\cdots u_{\MJ_k}$ with ${\MJ_k} \leq \cdots
\leq \MJ_1$ be two monomials in the minimal
generating set of $H_r(P)^{k}$ such that $\deg_{x_{m'{\ell}'}}u=\deg_{x_{m'{\ell}'}}v$ for all $x_{m'{\ell}'}>x_{m{\ell}}$ and $\deg_{x_{m{\ell}}}u>\deg_{x_{m{\ell}}}v$.

\medskip

We claim:
\begin{eqnarray}
\label{if}
I_{sm'}=J_{sm'} \quad \text{for all} \quad m'<m \quad\text{and all}\quad s.
\end{eqnarray}
We prove the claim by induction on $m'$. Let $m'=1$.  We have to show that $I_{s1}=J_{s1}$ for all $s$. Let $p_\ell\in P$. Then $p_\ell\in I_{s1}$ if and only if $\deg_{x_{1{\ell}}}u\geq s$. Similarly $p_\ell\in J_{s1}$ if and only if $\deg_{x_{1{\ell}}}v\geq s$. Since $\deg_{x_{1{\ell}}}u=\deg_{x_{1{\ell}}}v$, the desired conclusion follows.

Now let $m'<m$ and assume that $I_{sm''}=J_{sm''}$ for all $m''<m'$ and all $s$. Again let  $p_\ell\in P$. Then $p_\ell\in I_{sm'}\setminus I_{s,m'-1}$ if and only if $\deg_{x_{m'\ell}}u\geq s$, and similarly  $p_\ell\in J_{sm'}\setminus J_{s,m'-1}$ if and only if $\deg_{x_{m'\ell}}v\geq s$. Since $\deg_{x_{m'{\ell}}}u=\deg_{x_{m'{\ell}}}v$, it follows that $I_{sm'}\setminus I_{s,m'-1}=J_{sm'}\setminus J_{s,m'-1}$. Our induction hypothesis implies that $I_{s,m'-1}= J_{s,m'-1}$. Thus the desired conclusion follows.

Next we claim:
\begin{eqnarray}
\label{next}
\text{for all $\ell'<\ell$ and all $s$}, \quad p_{\ell'}\in I_{sm} \Leftrightarrow  p_{\ell'}\in J_{sm}.
\end{eqnarray}
As in the proof of claim (\ref{if}) we see that $p_{\ell'}\in I_{sm}\setminus  I_{s,m-1}$ if $\deg_{x_{m\ell'}}u\geq s$, and
$p_{\ell'}\in J_{sm}\setminus  J_{s,m-1}$ if $\deg_{x_{m\ell'}}v\geq s$. Hence since $\deg_{x_{m\ell'}}u=\deg_{x_{m\ell'}}v$ for $\ell'<\ell$, we conclude that $p_{\ell'}\in I_{sm}\setminus  I_{s,m-1}$ if and only if $p_{\ell'}\in J_{sm}\setminus  J_{s,m-1}$. However by claim (\ref{if}) we have that $I_{s,m-1}=J_{s,m-1}$. Thus the result follows.

\medskip
Since $\deg_{x_{m,\ell}}u> \deg_{x_{m,{\ell}}}v$, there exists some $t$ such that  $p_{\ell}\in I_{t,m}\setminus I_{t,m-1}$ and $p_{\ell}\not\in J_{t,m}\setminus J_{t,m-1}$. By (\ref{if}), we have  $I_{t,m-1}=J_{t,m-1}$. Therefore  $p_{\ell}\in J_{t,m'}\setminus J_{t,m'-1}$ for some $m'>m$.

Now, we consider the multichain $\mathcal{L}:\ L_1\subseteq L_2\subseteq \cdots \subseteq L_r=P$ of subsets of $P$, where
\[
L_{1}=J_{t,1},\ldots,L_{m-1}=J_{t,m-1}\quad
\]
\[
 L_{m}=J_{t,m}\cup\{p_{\ell}\},\ldots,L_{m'-1}=J_{t,m'-1}\cup\{p_{\ell}\},
\]
\[
L_{m'}=J_{t,m'},\ldots,L_{r}=J_{t,r}
\]
Observe that $\mathcal{L}$ is indeed a multichain of poset ideals of $P$. We have already $L_r=J_{t,r}=P$. Therefore, it is enough to show that $L_j$ is a poset ideal for each $j$. If $j<m$ or $j\geq m'$, then $L_{j}=J_{t,j}$ is a poset ideal. Assume that $m\leq j< m'$, then $L_{j}=J_{t,j}\cup\{p_{\ell}\}$.
Since $I_{t,m}$ is a poset ideal containing $p_{\ell}$, for any element $p_{\ell'}<p_{\ell}$ we have  $p_{\ell'}\in I_{t,m}$,  and so by (\ref{next}) we see that  $p_{\ell'}\in J_{t,m}$. Hence the monomial $u_{\MJ_1}\cdots u_{\MJ_{t-1}}u_{L}u_{\MJ_{t+1}}\cdots u_{\MJ_{k}}=x_{m,\ell}(v/x_{m',\ell})$ is a monomial in $H_r(P)^k$ which fulfills the condition of weakly polymatroidal ideals.
\end{proof}

We shall use Theorem~\ref{power} to show that the ideals $H_{r,s}(P)$ are weakly polymatroidal. For the proof of this fact we need the following simple result. Let $I$ be a monomial ideal generated in one degree. The {\em squarefree part} of $I$ is the ideal generated by all squarefree generators of $I$.

\begin{Lemma}\label{squarefree}
The squarefree part of every weakly polymatroidal ideal is again  weakly polymatroidal.
\end{Lemma}

\begin{proof}
Let $I$ be a weakly polymatroidal ideal in $K[x_1,\ldots,x_n]$. Let $u=x_{i_1}x_{i_2}\cdots x_{i_d}$ and $v=x_{j_1}x_{j_2}\cdots x_{j_{d}}$ be two monomials in the minimal generating set of the squarefree part of $I$ with $i_1=j_1,\ldots,i_{t-1}=j_{t-1}$ and $i_t<j_t$.  Since $I$ is weakly polymatroidal, there exists some $\ell\geq t$ such that  $w=x_{i_t}(v/x_{j_{\ell}})$ is in $I$. Since $w$ is again a squarefree monomial, it follows that the squarefree part of $I$ fulfills the condition of weakly polymatroidal ideals.
\end{proof}

Now we are ready to show

\begin{Theorem}
\label{irscm}
The ideal
$H_{r,s}(P)$ is weakly polymatroidal.
In particular, $I_{r,s}(P)$ is a Cohen--Macaulay ideal.
\end{Theorem}

\begin{proof} In Theorem~\ref{FatemehViviana} we observed that $H_{r,s}(P)=H_r(P)^{\langle r-s+1\rangle}$. Therefore  Theorem~\ref{power} and Lemma~\ref{squarefree} imply that $H_{r,s}(P)$ is weakly polymatroidal.

By \cite[Theorem~1.4]{Hibi2} weakly polymatroidal ideals have linear quotients. Thus, since $H_{r,s}(P)$ is weakly polymatroidal, it
follows from \cite[Proposition~8.2.5]{book} that  $I_{r,s}(P)$ is Cohen--Macaulay.
\end{proof}

\begin{Corollary}
\label{ircm}
The ring $S/I_r(P)$ is Cohen--Macaulay of dimension $n(r-1)$, and the following conditions are equivalent:
\begin{enumerate}
\item[(a)] $S/I_r(P)$ is Gorenstein;
\item[(b)] $S/I_r(P)$ is a complete intersection;
\item[(c)] all elements of $P$ are incomparable.
\end{enumerate}
\end{Corollary}

\begin{proof}
The degrees of the minimal generators of $H_r(P)$ correspond to the heights of the minimal prime ideals of $I_r(P)$. Since all generators of  $H_r(P)$ are of degree $n$ it follows that $\height I_r(P) =n$. Thus $\dim S/I_r(P)= nr-n=n(r-1)$.

Consider the sequence $\bf{s}$ of elements $x_{1j}-x_{ij}$ with $i=2,\cdots, r$ and $j=1,\ldots,n$. Then
$S/I_r(P)/({\bf s})S/I_r(P)\iso K[x_{11},\ldots, x_{1n}]/J$ where $J$ is generated by the monomials $x_{1i_1}\cdots x_{1i_r}$ with $p_{i_1}\leq p_{i_2}\leq \cdots \leq p_{i_r}$. Thus $J$ contains the elements $x_{11}^r,\ldots,x_{1n}^r$. In particular, $\dim  S/I_r(P)/({\bf s})S/I_r(P)=0$, which implies that ${\bf s}$ is a regular sequence, since the length of   ${\bf s}$ is $n(r-1)$.  It follows that $S/I_r(P)$ is Gorenstein if and only if $K[x_{11},\ldots, x_{1n}]/J$ is Gorenstein. Since $J$ is a monomial ideal and since $\dim K[x_{11},\ldots, x_{1n}]/J=0$, this is the case if and only if $J$ is generated by pure powers of the variables.  This happens if and only if $I_r(P)$ is generated by the monomials $x_{1i}\cdots x_{ri}$ with $i=1,\ldots,n$. This yields the desired conclusions.
\end{proof}

Each of the ideals $I_r(P)$ may be considered as the edge ideal of an uniform Cohen-Macaulay admissible clutter (see \cite{Villa2}, \cite{Villa1}). However, there exist uniform Cohen-Macaulay admissible clutters which do not arise from a poset. Such an example is given in \cite[Example 3.4]{Villa1}. Namely, the ideal $I=(x_1y_1z_1,x_2y_2z_2,x_1y_1z_2)\subset K[x_1,y_1,z_1,x_2,y_2,z_2]$ is the edge ideal of an uniform Cohen-Macaulay admissible clutter.  If $I$ came from a poset $P$, then  $P$ would have $2$ elements since $\height(I)=2,$ thus $I=I_3(P).$ But in this case $I$ cannot have three minimal monomial generators.

\section{The resolution of $H_r(P)$}

We recall that, for a monomial ideal $I\subset R=K[x_1,\ldots,x_n]$ which has linear quotients with respect to the order $u_1,\ldots,u_m$
of its minimal generators, we denote by $\set(u_i)$, for $i\geq 1$, the set of all the variables which generate the quotient ideal
$(u_1,\ldots,u_{i-1}): u_{i}.$ By \cite[Lemma 1.5]{HT}, the symbols $f(\sigma; u),$ $\sigma\subset \set(u),$ $|\sigma|=i-1,$ $u\in
G(I),$ form a homogeneous basis of the $i$th module in the minimal resolution of $R/I.$ Therefore, the computation of the sets of $I$ allows
the computation of the Betti numbers of $R/I$.

We now compute the sets associated with the minimal generators of $H_r(P).$

\begin{Lemma}\label{sets}
Let $\Ic \in \Ic_r(P)$ be a multichain of ideals in $P,$ $\Ic: I_1\subseteq \ldots \subseteq I_{r-1}\subseteq I_r=P$, and
$u_{\Ic}$ the corresponding generator of $H_r(P).$ Then
\[
\set(u_{\Ic})=\Union_{m=1}^{r-1}\{x_{mj}: p_j\in \Min(P\setminus I_m)\},
\]
where, for a poset ideal $I\subset P,$ $\Min(P\setminus I)$ is the set of all the minimal elements in $P\setminus I.$
In particular, we have
\[
\reg S/I_r(P)=\projdim H_r(P) = \max{\sum_{m=1}^{r-1}|\Min(P\setminus I_m)|},
\]
where the maximum is taken over all the multichains of poset ideals $\Ic: I_1\subseteq \ldots \subseteq I_{r-1}\subseteq I_r=P$.
\end{Lemma}

\begin{proof}
Let $x_{mj}$ with $p_j\in \Min(P\setminus I_m)$. It follows that there exists $\ell > j$ such that $p_j\in I_\ell.$
 Let $t=\min\{\ell : p_j\in I_\ell\},$ that is,  $p_j\in I_t\setminus I_{t-1}.$ Then
\[
\Ic^{\prime}: I_1^{\prime}=I_1\subseteq \ldots \subseteq I_{m-1}^{\prime}=I_{m-1}\subseteq I_m^{\prime}=I_m\union \{p_j\}\subseteq\cdots\subseteq
I_{t-1}^{\prime}=I_{t-1}\union \{p_j\}\subseteq
\]
\[
\subseteq I_t^{\prime}=I_t\subseteq \ldots \subseteq I_r^{\prime}=I_r=P.
\]
is a multichain of poset ideals, $u_{\Ic^{\prime}}=x_{mj}(u_\Ic/x_{tj}),$ and $u_{\Ic^{\prime}}>_{\lex} u_{\Ic}.$ Therefore,
\[
x_{mj}=u_{\Ic^{\prime}}/ \gcd(u_{\Ic^{\prime}},u_\Ic)\in \set(u_{\Ic}).
\]

For the other inclusion, let $u_{\Lc}>_{\lex} u_\Ic.$ Thus there exist $m$ and $j$ such that
$\deg_{x_{m^\prime\j^\prime}}u_\Lc=\deg_{x_{m^\prime\j^\prime}}u_\Ic$ for all $x_{m^\prime\j^\prime} > x_{mj}$ and
$\deg_{x_{mj}}(u_\Lc)> \deg_{x_{mj}}(u_\Ic).$ Thus $p_j\not\in I_m$. As in the first part of the proof we define the monomial
$u_{\Ic^\prime}$ and  we get that  $u_{\Ic^\prime}=x_{mj}u_\Ic / x_{tj} >_{\lex} u_{\Ic}.$ Obviously $x_{mj}$ divides
$u_{\Lc}/\gcd(u_\Lc,  u_\Ic)$, which ends the proof.

The formula for the projective dimension is an immediate consequence \cite[Lemma 1.5]{HT}, while the equality with the regularity is implied by
a result of Terai (see \cite[Prop.~8.1.10]{book}).
\end{proof}

In a similar way one can prove the following slightly more general

\begin{Lemma}
\label{moregeneral}
Let $\MI_1,\ldots,\MI_k \in \Ic_r(P)$ be  multichains of ideals in $P$ with $\MI_l: I_{l1}\subseteq \ldots \subseteq I_{lr-1}\subseteq I_{lr}=P$, and
$u_{\MI_1}\cdots u_{\MI_k}$ the corresponding generator of $H_r(P)^k.$ Then
\[
\set(u_{\MI_1}\cdots u_{\MI_k})=\Union_{l=1}^k\Union_{m=1}^{r-1}\{x_{mj}: p_j\in \Min(P\setminus I_{lm})\}.
\]
\end{Lemma}

Recall that an antichain of $P$ is a subset $A$ of $P$ which any two of whose elements are
incomparable in $P$. By using this concept we get the following interpretation of the regularity of $I_r(P)$.
\begin{Corollary}
\label{H_r(P)}
We have \[
\reg S/I_r(P)= \projdim H_r(P) = (r-1)s,
\] where $s$ is the maximal  cardinality of an  antichain of $P$.
\end{Corollary}
\begin{proof}
Let $A\subseteq P$ be an antichain with $|A|=s$,  and let $B$ be the following ideal of $P$:
\[
B=\{q\in P\;\: q<p\quad \text{for some $p\in A$}\}.
\]
Now, consider the following multichain of poset ideals
\[
\MI:\ I_{1}=B = \cdots=I_{r-1}=B\subset I_{r}=P.
\]
Then
\[
\sum_{m=1}^{r-1}|\Min(P\setminus I_m)|=(r-1)s,
\]
and obviously  the maximum of the numbers $\sum_{m=1}^{r-1}|\Min(P\setminus I'_m)|$ taken over all multichains $\MI'$ of $\MI(P)$ cannot be bigger than that for $\MI$. The desired conclusion follows from Lemma~\ref{sets}.
\end{proof}

\medskip
In order to determine the resolution of $H_r(P)$ we apply the method developed in \cite{HT}. We first recall the needed tools. For
a  monomial ideal $I\subset R=K[x_1,\ldots,x_n]$ which has linear quotients with respect to the order $u_1,\ldots,u_m$ of its minimal
monomial generators, one defines its decomposition function $g: M(I)\to G(I)$ by $g(u)=u_j$ if $j$ is the smallest number such that
$u\in (u_1,\ldots,u_{j-1}):u_j.$ Here $M(I)$ denotes the set of all monomials of the ideal $I.$ The function $g$ is regular if
$\set(g(x_su))\subset \set(g(u))$ for all $s\in \set(u)$ and $u\in G(I).$

In order to show that the decomposition function associated with $H_r(P)$ is regular, we fix a notation. For $u_{\Ic}\in G(H_r(P))$ and
$x_{mj}\in \set(u_\Ic)$ we denote by $u_{\Ic^{\prime}}$ the generator of $H_r(P)$ corresponding to the following multichain of poset ideals
\[
\Ic^{\prime}: I_1^{\prime}=I_1\subseteq \ldots \subseteq I_{m-1}^{\prime}=I_{m-1}\subseteq I_m^{\prime}=I_m\union \{p_j\}\subseteq\cdots\subseteq
I_{t-1}^{\prime}=I_{t-1}\union \{p_j\}\subseteq
\]
\[
\subseteq I_t^{\prime}=I_t\subseteq \ldots \subseteq I_r^{\prime}=I_r=P.
\]
where, as before, $t=\min\{\ell: p_j\in I_\ell\}$. We have seen in the proof of Lemma \ref{sets} that
$u_{\Ic^{\prime}}=x_{mj}(u_\Ic/x_{tj})$ and $u_{\Ic^{\prime}}>_{\lex} u_{\Ic}.$

\begin{Lemma}\label{decomp1}
Let $g$ be the decomposition function of $H_r(P),$ $u_\Ic$ a minimal generator of $H_r(P)$ and $x_{mj}\in \set(u_\Ic).$ Then
\[
g(x_{mj}u_\Ic)=u_{\Ic^\prime}.
\]
\end{Lemma}

\begin{proof}
Let $g(x_{mj}u_\Ic)=u_\Lc$ for some multichain of poset ideals $\Lc: L_1\subseteq \ldots\subseteq L_{r-1}\subseteq L_r=P.$ This implies that
$x_{mj}u_\Ic=u_\Lc x_{\nu q}$ for some variable $x_{\nu q}.$ Let us suppose that $x_{\nu q}\neq x_{tj}.$ Since $x_{mj} | u_{\ML}$ it
follows that $p_j\in L_m\setminus L_{m-1}$. On the other hand, since
$x_{tj} | u_\MI$ and $x_{\nu q}\neq x_{tj},$ it follows that $x_{tj} | u_\ML,$ thus $p_j\in L_t\setminus L_{t-1}.$ Therefore, we must
have  $t= m,$ which is impossible.
 \end{proof}

\begin{Corollary}\label{decomp2}
\[
\set(x_{mj}u_\Ic)\subset \set(u_\Ic).
\]
\end{Corollary}

\begin{proof}
By the definition of the multichain $\Ic^{\prime},$ we have that $I_\mu^\prime=I_\mu$ for all $\mu < m$ or $\mu \geq t.$ Let $m\leq \mu < t.$
Then $\Min (P\setminus I_\mu^{\prime})=\Min (P\setminus I_\mu)\setminus \{p_j\}\subset \Min(P\setminus I_\mu)$. By applying  Lemma \ref{sets},
 we get the inclusion.
\end{proof}

\cite[Theorem~1.12]{HT} gives the minimal resolution of  monomial ideals with linear quotients which admit a regular decomposition function. Applying this theorem to our situation we obtain

\begin{Theorem}
\label{resolution} Let $\mathbb{F}_{\bullet}$ be the graded minimal free resolution of $S/H_r(P)$. Then the symbols $f(\sigma; u_\Ic),$ $\sigma\subset
\set(u_\Ic),$ $|\sigma|=i-1,$ $u_\Ic\in  G(H_r(P)),$ form a homogeneous basis of $F_i$ for $i\geq 0.$ The chain map of $\mathbb{F}_{\bullet}$ is given by
\[
\partial(f(\sigma; u_\Ic))=\sum_{x_{mj}\in \sigma}(-1)^{\alpha(\sigma,x_{mj})}(x_{tj}f(\sigma\setminus\{x_{mj}\};u_{\Ic^\prime})-
x_{mj}f(\sigma\setminus\{x_{mj}\};u_{\Ic}))
\]
if $\sigma\neq\emptyset$, and
\[\partial(f(\emptyset; u_\Ic)) = u_\Ic\] otherwise.
Here $\alpha(\sigma, x_{mj})=|\{x_{m^\prime \ell^\prime}\in\sigma\ : x_{m^\prime \ell^\prime}> x_{mj} \}|$.
\end{Theorem}

\section{Generalized Hibi rings}

In this section we introduce a class of $K$-algebras which can be identified with classical Hibi rings. Let  $R_r(P)$ be the toric ring  generated over $K$ by all the monomials $u_\MI$ with $\MI\in\mathcal{I}_r(P)$, and let $T$ be the polynomial ring over $K$ in the indeterminates $y_{\MI}$ with $\MI\in \mathcal{I}_r(P)$.  Furthermore let $\varphi\: T\to R$ be the surjective $K$-algebra homomorphism with $\varphi(y_{\MI})=u_{\MI}$ for all $\MI \in \mathcal{I}_r(P)$.
We choose  a total order on the variables $y_{\MI}$ with the property that $\MI<  \MI'$ implies that $y_{\MI}> y_{\MI'}$.

\begin{Theorem}
\label{generalized}
The set $\mathcal G$ of elements
\[
y_{\MI}y_{\MI'}-y_{\MI\union \MI'}y_{\MI\sect \MI'}\in T\quad  \text{with}\quad \MI,\MI'\in \mathcal{I}_r(P)\quad \text{incomparable,}
\]
is  a  reduced Gr\"obner basis of the ideal $L_r=\Ker\varphi$ with respect to the reverse lexicographic order induced by the given order of the variables $y_{\MI}$.
\end{Theorem}

\begin{proof} The equation (\ref{unionsect}) shows that $y_{\MI}y_{{\MI}'}-y_{{\MI}\union {\MI}'}y_{{\MI}\sect {\MI}'}$ is indeed an element of $L_r$.
Now let $\prod_{s=1}^ty_{\MI_s}-\prod_{s=1}^ty_{\MI'_s}$ be a primitive  binomial  in $L_r$ with initial monomial $\prod_{s=1}^ty_{\MI_s}$. We are going to show that there are two indices $k$ and $\ell$  such that $\MI_k$ and $\MI_{\ell}$ are incomparable multichains of ideal, and that $y_{\MI_k}y_{\MI_{\ell}}$ is the leading monomial of $y_{\MI_k}y_{\MI_{\ell}}-y_{\MI_k\union \MI_{\ell}}y_{\MI_k\sect \MI_{\ell}}$. This will then show that $\mathcal G$ is  Gr\"obner basis of $L_r$. It is obvious that $\mathcal G$  is actually reduced.

Suppose to the contrary that $\MI_1\leq \MI_2\leq \ldots \leq \MI_t$. We will show that $\MI'_s< \MI_t$ for all $s$. Indeed, since $\prod_{s=1}^ty_{\MI_s}-\prod_{s=1}^ty_{\MI'_s}\in L_r$ we see  that $\prod_{s=1}^tu_{\MI_s}=\prod_{s=1}^tu_{\MI'_s}$. It follows that
\[
\prod_{s=1}^t(\prod_{k=1}^{\ell}x_{k I_{sk}\setminus I_{sk-1}})=
\prod_{s=1}^t(\prod_{k=1}^{\ell}x_{k I'_{sk}\setminus I'_{sk-1}})\quad \text{for all $\ell=1,\cdots,r$}.
\]
Here $\MI_s$ is the multichain of ideals  $I_{s1}\subseteq I_{s2}\subseteq \cdots \subseteq I_{sr}=P$,  and  $\MI'_s$ the multichain of ideals $I'_{s1}\subseteq I'_{s2}\subseteq \cdots \subseteq I'_{sr}=P$.

Now for all $j$ and $k$ we apply the substitution $x_{kj}\mapsto x_j$, and obtain
\[
\prod_{s=1}^tx_{I_{s{\ell}}}=\prod_{s=1}^tx_{I'_{s{\ell}}}, \quad {\ell}=1,\ldots,r,
\]
where $x_J=\prod_{j\in J}x_j$ for $J\subset [n]$.

Since $\MI_1\leq \MI_2\leq \ldots \leq \MI_t$, it follows that $\supp(\prod_{s=1}^tx_{I_{s{\ell}}})=I_{t{\ell}}$. Thus the equation $\prod_{s=1}^tx_{I_{s{\ell}}}=\prod_{s=1}^tx_{I'_{s{\ell}}}$ implies that $x_{I'_{s{\ell}}}|x_{I_{t{\ell}}}$ for all ${\ell}$ and all $s$. It follows that $\MI'_s\leq \MI_t$. We cannot have equality, since   $\prod_{s=1}^ty_{\MI_s}-\prod_{s=1}^ty_{\MI'_s}$ is a primitive  binomial.  This contradicts the fact that  $\prod_{s=1}^ty_{\MI_s}$ is the initial  monomial of  $\prod_{s=1}^ty_{\MI_s}-\prod_{s=1}^ty_{\MI'_s}$.

Finally, $y_{\MI_k}y_{\MI_{\ell}}$ is the leading monomial of $y_{\MI_k}y_{\MI_{\ell}}-y_{\MI_k\union \MI_{\ell}}y_{\MI_k\sect \MI_{\ell}}$ thanks to the monomial order on $T$.
\end{proof}

\begin{Corollary}
\label{normalcm}
For any poset $P$ and all integers $r\geq 1$, the toric ring $R_r(P)$ is a normal Cohen--Macaulay domain.
\end{Corollary}

\begin{proof}
Since the defining ideal of $R_r(P)$ has a squarefree initial ideal, it follows from a result of Sturmfels \cite[Corollary 8.8]{St} that $R_r(P)$ is normal, and a result of Hochster  \cite[Theorem 1]{Ho} that $R_r(P)$ is Cohen--Macaulay.
\end{proof}

Our next goal is to find out which of the toric rings $R_r(P)$ is Gorenstein. The answer to this will be a consequence of the next theorem where it will be shown that $R_r(P)$ can be interpreted as Hibi ring of a suitable poset. To be precise, let $Q_r=[r]$ with usual order. Recall that the direct product $P\times Q$ of two poset $P$ and $Q$ is poset on product of the underlying sets of $P$,  $Q$ with partial order given by
\[
(p,q)\leq (p',q')\quad  \Leftrightarrow \quad p\leq p' \quad\text{and}\quad q\leq q'.
\]
\begin{Theorem}
\label{hibiring}
Let $P$ be any finite poset. Then $R_r(P)\iso R_2(P\times Q_{r-1})$ for all $r\geq 2$.
\end{Theorem}

\begin{proof}
We first notice that the partially ordered set $\MI_r(P)$ is a distributive lattice with meet and join defined as intersection  and union of multichains of ideals. In Theorem~\ref{generalized} we have seen that the defining relations of $R_r(P)$ are just the Hibi relations of the distributive lattice $\MI_r(P)$. In particular, it follows that $R_r(P)$ is the Hibi ring of $\MI_r(P)$. Let $P'$ be the subposet of join irreducible elements in $\MI_r(P)$. Then we obtain that $R_r(P)\iso R_2(P')$. Thus it remains to be shown that $P'\iso P\times Q_{r-1}$. For this purpose we have to identify the join irreducible elements in  $\MI_r(P)$.

Let $\MI\: \emptyset \subseteq \cdots \subseteq \emptyset \subset I_k\subseteq \cdots \subseteq I_r=P$ be an element of  $\MI_r(P)$. We  claim that
$\MI$ is join irreducible if and only if $I_k$ is a join irreducible element of the ideal lattice of $P$ and $I_k=I_{k+1}=\cdots =I_{r-1}$. Indeed, suppose that $I_k=J\union J'$ where $J$ and $J'$ are ideals in $P$ properly contained in $I_k$. Then $\MI= \MJ\union \MJ'$ where
\[
\MJ\: \emptyset \subseteq \cdots \subseteq \emptyset \subset J\subseteq I_{k+1}\subseteq \cdots \subseteq I_r
\]
and
\[
\MJ'\: \emptyset \subseteq \cdots \subseteq \emptyset \subset J'\subseteq I_{k+1}\subseteq \cdots \subseteq I_r,
\]
a contradiction.

Suppose now that there exists an integer $s$ with $k\leq s<r-1$ such that $I_j=I_k$ for $j=k,\ldots,s$  and $I_k\subset I_{s+1}$. Then $\MI= \MJ\union \MJ'$ where
\[
\MJ\: \emptyset \subseteq \cdots \subseteq \emptyset \subset I_k\subseteq I_{k}\subseteq \cdots \subseteq I_k\subseteq I_r
\]
and
\[
\MJ'\: \emptyset \subseteq \cdots \subseteq \emptyset \subset I_s \subseteq I_{s+1}\subseteq \cdots \subseteq I_r,
\]
a contradiction.

Thus we have shown the "only if" part of the claim. The "if" part is obvious.

Let $\MI: \emptyset\subseteq \cdots\subseteq \emptyset \subset I=I=\cdots=I\subset P$ with $k$ copies of $I$ where $I$ is a join irreducible element of $\MI(P).$ Then $I$ is a principal ideal in $\MI(P),$ hence there exists a unique element $p\in I$ such that $I=\{a\in P : a\leq p\}.$ 

Finally we define the  poset isomorphism between the poset of join irreducible elements of $\MI_r(P)$ and $P\times Q_{r-1}$ as follows. To $\MI: \emptyset\subseteq \cdots \subseteq\emptyset \subset I=I=\cdots=I\subset P$ with $k$ copies of $I$ we assign $(p,k)\in P\times Q_{r-1}$.\end{proof}

\begin{Corollary}
\label{dim}
Let $P$ be poset of cardinality $n$. Then $\dim R_r(P)=n(r-1)+1$.
\end{Corollary}

\begin{proof} It is known \cite{Hibi1} and easy to see that the classical Hibi ring on a poset of cardinality $m$ has Krull dimension $m+1$. Since  $R_r(P)\iso R_{2}(P\times Q_{r-1})$ and since $|P\times Q_{r-1}|=n(r-1)$ the assertion follows.
\end{proof}

A poset $P$ is called {\em pure}, if all maximal chains have the same length.

\begin{Corollary}
\label{hibiringgor}
Let $P$ be a finite poset. Then the following conditions are equivalent:
\begin{enumerate}
\item[(a)] $R_r(P)$ is Gorenstein;
\item[(b)] $R_2(P)$ is Gorenstein;
\item[(c)] $P$ is pure.
\end{enumerate}
\end{Corollary}

\begin{proof}  A well-known theorem of Hibi \cite[Corollary 3.d]{Hibi1} says that $R_2(P')$ is Gorenstein if and only if $P'$ is pure. Since $P$ is pure if and only if $P\times Q_{r-1}$ is pure, it follows that all the statements are equivalent.
\end{proof}

\section{The algebra $R_{r,s}(P)$}

Now we want to study the algebra $R_{r,s}(P)$ and show that it has a quadratic Gr\"obner basis with squarefree initial ideal. Let $G_{r,s}$ be the minimal set of monomial generators of $H_{r,s}(P)$. Then the elements of $G_{r,s}$ generate $R_{r,s}(P)$. Let $k=r-s+1$, then $G_{r,s}$ consist of all squarefree monomials of the form $u_{\MI_1}u_{\MI_2}\cdots u_{\MI_k}$ with ${\MI_1}> \MI_2>\cdots>  \MI_k$ and $\MI_j\in \MI_r(P)$. Corresponding to each such monomial we introduce the variable $y_{\MI_1,\MI_2,\ldots, \MI_k}$, and let $T$ be the polynomial ring over $K$ in this set of variables. Let $L_{r,s}=\Ker \varphi$ where $\varphi\: T\to R_{r,s}(P)$ is the $K$-algebra homomorphism with $y_{\MI_1,\MI_2,\ldots, \MI_k}\mapsto u_{\MI_1}u_{\MI_2}\cdots u_{\MI_k}$.

Defining in a similar  way the order of the variables $y_{\MI_1,\MI_2,\ldots, \MI_k}$ as we did it in the case of $R_r(P)$,  we choose a total order on the variables with the property that
\[
 y_{\MI_1,\MI_2,\ldots, \MI_k}>  y_{\MJ_1,\MJ_2,\ldots, \MJ_k}\quad \text{if}\quad (\MI_1,\MI_2,\ldots, \MI_k)<(\MJ_1,\MI_2,\ldots, \MJ_k),
\]
where by definition $(\MI_1,\MI_2,\ldots, \MI_k)\leq (\MJ_1,\MJ_2,\ldots, \MJ_k)$, if $\MI_l\leq \MJ_l$ for $l=1,\ldots,k$.

In analogy to Theorem~\ref{generalized} one would expect that $L_{r,s}$ has a quadratic Gr\"obner basis with respect to the reverse lexicographic order induced by the above order of the variables. This is however not the case. To see this we choose for $P$ the poset consisting of only one element $p$. Then $H_{r,s}(P)=I_{r,s}$ where $I_{r,s}$ denotes the squarefree Veronese ideal of degree $s$ in $r$ variables, that is, the ideal generated by all squarefree monomials of degree $s$  in $T=K[x_1,\ldots,x_r]$. In this particular  case the above defined order of the variables is the  $y_u > y_v$ if $u>_{\lex} v$.

\begin{Theorem}
\label{veroneseviv}
Let $\mathcal{G}$ be the reduced Gr\"obner basis of $L_{r,s}=\Ker\varphi$ with respect to the reverse lexicographical order  on $T$ induced by the above order of the variables, where $\varphi\: T\to R_{r,s}(P)$ is the $K$ algebra homomorphism with $y_u\mapsto u$. Then the
binomials of $\mathcal{G}$ have squarefree initial monomials and are generated in degree at most $3.$
\end{Theorem}

\begin{proof}
Let $g=y_{u_1}\cdots y_{u_q}-y_{v_1}\cdots y_{v_q}\in \mathcal{G}$ with $u_1\geq_{\lex}\ldots \geq_{\lex}u_q $ and
$v_1\geq_{\lex}\ldots \geq_{\lex}v_q.$ Let $\ini_{<}(g)=y_{u_1}\cdots y_{u_q}.$ Then, since $g$ is a primitive binomial, we have
$y_{u_q}>y_{v_q},$ that is, $u_q>_{\lex}v_q.$

Obviously, there are at least two different variables in the support of $\ini_<(g),$ that is, $q\geq 2$ and $u_1>_{\lex}u_q.$

In the first place we assume that there exist $1\leq a < b\leq q$ with $u_a\neq u_b$ and $i\in \supp(u_a)\setminus \supp(u_b),$
$j\in \supp(u_b)\setminus \supp(u_a),$ such that $i> j.$ Let $u_a^{\prime}=x_ju_a/ x_i$ and $u_b^\prime=x_iu_b/x_j.$ Then we have
$u_au_b=u_a^{\prime}u_b^{\prime},$ that is $h=y_{u_a}y_{u_b}-y_{u_a^{\prime}}y_{u_b^\prime}\in L_{r,s}(P)$, and $u_b>_{\lex}u_b^\prime,$ whence
$\ini_<(h)=y_{u_a}y_{u_b}.$ Since $\mathcal{G}$ is a reduced Gr\"obner basis, we must have $\ini_<(g)=\ini_<(h),$ thus $\ini_<(g)$ is
a squarefree monomial of degree $2$.

Now we assume that for all $1 \leq a < b \leq q$ with $u_a\neq u_b$ we have
\begin{equation}\label{eq1}
\max(\supp(u_a)\setminus \supp(u_b))< \min(\supp(u_b)\setminus\supp(u_a)).
\end{equation}
Since $u_q>_{\lex}v_q$ we also have
\begin{equation}\label{eq2}
\ell=\min(\supp(u_q)\setminus\supp(v_q))< h=\min(\supp(v_q)\setminus\supp(u_q)).
\end{equation}

In particular, we get
\[
x_{\ell}| u_q, x_\ell\not| v_q, x_h | v_q, \text{ and } x_h\not|u_q.
\]
We then obtain that there exists $1\leq a,b\leq q-1$ such that $x_\ell\not| u_a$  and $x_h | u_b$. Note that $u_a\neq u_b.$ Indeed, if
$u_a=u_b$ we get $h\in \supp(u_a)\setminus \supp(u_q)$ and $\ell\in \supp(u_q)\setminus \supp(u_a).$ Since $a <q,$ by using (\ref{eq1}), we obtain $h < \ell$ which is in contradiction to (\ref{eq2}). Therefore, the monomials $u_a,u_b,$ and $u_q$ are distinct.
Let $j\in \supp(u_a)\setminus\supp(u_b)$ and consider the monomials
\[
u_a^\prime=x_\ell u_a/ x_j, u_b^\prime=x_j u_b/ x_h,  \text{ and } u_q^\prime=x_h u_q/ x_\ell.
\]
It follows that $u_au_bu_q=u_a^\prime u_b^\prime u_q^\prime$, thus
$h=y_{u_a}y_{u_b}y_{u_q}-y_{u_a^\prime}y_{u_b^\prime}y_{u_q^\prime}\in \mathcal{G}$, and $u_q>_{\lex} u_q^\prime$, hence
$\ini_<(h)=y_{u_a}y_{u_b}y_{u_q}.$ Since $\mathcal{G}$ is reduced, it follows that $\ini_<(g)$ is a squarefree monomial of degree at most $3.$
\end{proof}

The following example shows that the reduced Gr\"obner basis $\mathcal{G}$ of $L_{r,s}$ does in general indeed contain monomials of degree 3. With CoCoA we compute $\mathcal{G}$ for the ring $R_{6,3}(P)$. The binomials of degree 3 in  $\mathcal{G}$ are the following:  $kps - lmt, ejs - fgt, bjp - cdt, drs - gmt, cqs - flt, ajp - cds,  bqr - ekt, aqr - eks,  ano - bkp, aio - bdr, ahi - bej, ahn - bcq$.

\medskip
\noindent
Here, we denoted for simplicity the  variables of $T$ by $a,b,c,\cdots,t$ and defined $T\to R_{6,3}(P)$ by mapping the variables in their natural (lexicographical order)  to the corresponding  monomial of $H_{6,3}(P)$ in the lexicographic order. As can been seen, there are 12  binomials of degree 3 in the reduced Gr\"obner basis of $L_{6,3}$.

If $R_{6,3}(P)$ would be isomorphic to Hibi ring, then it would  have to have a quadratic Gr\"obner basis with respect to the reverse lexicographic order induced by some order of the variables. This is not the case, at least for the natural order of the variables $u$, as we have seen above, and very likely for any  other order of the variables. Unfortunately this is not so easy to check because there exist  quite a lot of different orders, even if one takes into account all the symmetries.

The following simple argument shows that the even smaller ring $R_{4,2}(P)$ with $P=\{p\}$  can not be a Hibi ring. It is generated over $K$ by the monomials
\[
x_1x_2,x_1x_3, x_1x_4, x_2x_3, x_2x_4, x_3x_4.
\]
Suppose $R_{4,2}(P)$ is the Hibi ring of a poset $P$. Then its ideal lattice $\MI(P)$ should have cardinality 6. The only posets with this property are
\begin{eqnarray*}
P_1&=&\{p_1,p_2,p_3,p_4,p_5\},\quad  p_1<p_2<p_3<p_4<p_5;\\
P_2&=&\{p_1,p_2,p_3, p_4\},\quad  p_1<p_2\quad \text{and}\quad p_2<p_3,p_4;\\
P_3&=&\{p_1,p_2,p_3, p_4\},\quad  p_1,p_2<p_3\quad  \text{and}\quad  p_3<p_4;\\
P_4&=&\{p_1,p_2,p_3\},\quad  p_1<p_2.
\end{eqnarray*}
For the posets $P_1,P_2,P_3$ the toric ideal of  the associated Hibi ring is generated by at most one binomial,  and for $P_4$ by three binomials,  while the defining ideal of $R_{4,2}(P)$ has two binomial generators.

\medskip
In order to obtain a squarefree Gr\"obner basis of the defining ideal of $R_{r,s}(P)$ we use a result of Sturmfels \cite[Theorem 14.2]{St} and show that  the  set $G_{r,s}$ is sortable.

Recall that a set $B$ of monomials which are of same degree $d$ in the polynomial ring $S=K[x_1,\ldots,x_n]$ is called {\em sortable}, if image of the map
\[
\sort\: B\times B\to S_d\times S_d
\]
is contained in $B\times B$. The map $\sort$ is defined as follows: let $u,v\in B$ where $uv=x_{i_1}x_{i_2}\cdots x_{i_{2d}}$ with $i_1\leq i_2\leq \cdots \leq i_{2d}$.  Then $\sort(u,v)=(u',v')$ where  $u'=x_{i_1}x_{i_3}\cdots x_{i_{2d-1}}$ and $v'=x_{i_2}x_{i_4}\cdots x_{i_{2d}}$. We call the pair of monomial {\em unsorted}, if $(u,v)\neq (u',v')$.

\begin{Theorem}[Sturmfels]
\label{sturmfels}
Let $R$ be the toric ring generated over $K$ by a sortable set $B$ of monomials. Then $R=K[y_u\: u\in B]/I$, and  the binomials
\[
y_uy_v-y_{u'}y_{v'}\quad \text{where  $(u, v)$ is unsorted and $(u',v')=\sort(u,v)$,}
\]
form a Gr\"obner basis of $I$.
\end{Theorem}

Theorem~\ref{sturmfels} will be used to prove

\begin{Theorem}
\label{fatemehsortable}
For all integers $1\leq s\leq r$ the set of monomials $G_{r,s}(P)$ is sortable. In particular, $R_{r,s}(P)$ has a quadratic Gr\"obner basis and is a normal Cohen--Macaulay domain.
\end{Theorem}

Before giving the proof of the theorem we need the following two lemmata.

\begin{Lemma}
\label{tonight}
Let $u_{\MI}=u_{\MI_1}\cdots u_{\MI_k}$ be an element in $G_{r,s}$ with
\[
{\MI_1}
> \cdots > \MI_k\quad \text{where}\quad
\MI_j \: I_{j1}\subset I_{j2}\subseteq \ldots \subseteq  I_{jr}\quad \text{for}\quad  j=1,\ldots,k,
\]
and let $\prod_{c\in C}x_{ct}|u_\MI$ for some $C\subseteq[r]$ and some $t\in[n]$. Then for any element $p_m<p_t$ of $P$, there exists a set $C'$ and a bijection $C\to C'$ with $c\mapsto c'\leq c$ such that
\[
\prod_{c'\in C'}x_{c'm}\ \text{ divides}\quad u_\MI.
\]
\end{Lemma}

\begin{proof}
First observe that $x_{ct}|u_\MI$ if and only if $p_t$ belongs to $I_{j c}\setminus I_{j , c-1}$ for some $j$. Note that $j$ is uniquely determined by $c$ and $t$, since $u_\MI$ is a squarefree monomial.
Then for pairwise disjoint indices $j_{1},\ldots,j_{{|C|}}$ we have $p_t\in I_{j_{i}, c_i}\setminus I_{j_{i}, c_i-1}$ for all $i$. Since $I_{{j_i}}$ is a poset ideal,
there exists some
$c_i'\leq c_i$ with
$p_m\in I_{j_{i}, c_i'}\setminus I_{j_{i}, c_i'-1}$. Then our first observation shows that $x_{c_i',m}|u_{\MI}$. Since $u_{\MI}$ is a squarefree monomial, $c_1',c_2',\ldots,c'_{|C|}$ are again pairwise disjoint indices and so $\prod_{c'\in C'}x_{c'm}$ divides $u_\MI$, where $C'=\{c_1',c_2',\ldots,c'_{|C|}\}$.
\end{proof}

\medskip

\begin{Lemma}
\label{terrible}
Let $u=x_{1,A_1}x_{2,A_2}\cdots x_{r,A_r}\in S_{kn}$ be a squarefree monomial satisfying the following condition $(*)$: for each $j\in [n]$ there exist exactly $k$ of the sets $A_i$, say $A_{i_1},\ldots,A_{i_k}$,  such that $j\in A_{i_l}$ for $l=1,\ldots,k$.

Then
$u=u_1u_2\cdots u_k$ where $u_i=x_{1,A_{i1}}x_{2,A_{i2}}\cdots x_{r,A_{ir}}$ such that
\begin{enumerate}
\item[(1)] for each $i=1,\ldots,k$ and for each $j\in[n]$ there exists a unique $l\in [r]$ such that  $j\in A_{il}$;
\item[(2)] $A_{i+1,j}\subseteq \Union_{l=1}^{j-1}A_{il}$ for all $i$  and $j$.
\end{enumerate}
\end{Lemma}

\begin{proof}
Let $A_{1i}= A_i\setminus \Union_{l=1}^{i-1}A_l$. Since $[n]$ is the disjoint union of the sets $A_{11},\ldots, A_{1r}$, condition (1) is satisfied for $i=1$. Let $v=u/u_1=x_{1,B_1}\cdots x_{r,B_r}$. Then $v$ is a squarefree monomial of degree $(k-1)n$, and  since condition (1) is satisfied for $i=1$ it follows from $(*)$   that for each $j\in [n]$ there exist exactly $k-1$ of the sets $B_i$, say $B_{i_1},\ldots,B_{i_{k-1}}$,  such that $j\in B_{i_l}$ for $l=1,\ldots,k-1$. By using induction on $k$ we may assume that $v=u_2u_3\cdots u_k$ with $u_i=x_{1,A_{i1}}x_{2,A_{i2}}\cdots x_{r,A_{ir}}$ such that the conditions (1) and (2) are satisfied for $i\geq 2$. Thus it remains to be shown that $A_{2j}\subseteq  \Union_{l=1}^{j-1}A_{1l}$ for all $j$. We actually show that $B_j\subset \Union_{l=1}^{j-1}A_{1l}$ for all $j$. Indeed,
\[
B_j=A_j\setminus A_{1,j-1}=
A_j\sect \Union_{l=1}^{j-1}A_l\subseteq \Union_{l=1}^{j-1}A_l= \Union_{l=1}^{j-1}A_{1l}.
\]
\end{proof}

First observe that with the similar notation as in the above lemma one has
\begin{eqnarray}
\label{divides}
 t\in A_{ls} \text{ for some $l$ } \Leftrightarrow x_{st} \text{ divides } u.
\end{eqnarray}

\begin{proof}[Proof of Theorem~\ref{fatemehsortable}]
By  Theorem~\ref{FatemehViviana} we have $H_{r,s}(P)=H_r(P)^{\langle k\rangle}$ where $k=r-s+1$. Consider the following order
\[
x_{11}>x_{21}>\cdots>x_{r1}>x_{12}>\cdots>x_{r2}>\cdots>x_{1n}>\cdots>x_{rn}
\]
of the variables of $S$.

Let $u_{\MI}=u_{\MI_1}\cdots
u_{\MI_k}$ and $u_{\MJ}=u_{\MJ_1}\cdots
u_{\MJ_k}$ be two elements in $G_{r,s}$ with
\[
{\MI_1}
> \cdots > \MI_k\quad \text{and}\quad  {\MJ_1}
 >\cdots > \MJ_k,
 \]
where
\[
\MI_j \: I_{j1}\subseteq I_{j2}\subseteq \ldots \subseteq  I_{jr}\quad  \text{and}\quad
\MJ_j \: J_{j1}\subseteq J_{j2}\subseteq \ldots \subseteq  J_{jr}\quad  \text{for}\quad  j=1,\ldots,k.
\]
\medskip
Let
\[
\sort(u_{\MI},u_{\MJ})=(u,v).
\]
We first notice that both monomials $u$ and $v$ are again squarefree. Indeed, since $u_{\MI}u_{\MJ}=uv$ it follows that  each variable $x_{ij}$ appears at most to the power 2 in $uv$. If this happens, then the sorting operator moves one $x_{ij}$ to $u$ and the other $x_{ij}$  to $v$.

Now we may decompose $u$ and $v$ as in Lemma~\ref{terrible}. Say,
\[
u=u_1u_2\cdots u_k\quad\text{with}\quad u_i=x_{1,A_{i1}}x_{2,A_{i2}}\cdots x_{r,A_{ir}},
\]
and
\[
v=v_1v_2\cdots v_k\quad\text{with}\quad v_i=x_{1,B_{i1}}x_{2,B_{i2}}\cdots x_{r,B_{ir}}.
\]
The proof of the theorem is completed once we have shown that for all $i$ and $j$ the sets
\[
\Union_{l=1}^jA_{il}\quad\text{and}\quad  \Union_{l=1}^jB_{il}
\]
are poset ideals in $P$.

\medskip
Let $t\in A_{ij}$ and suppose that $p_m<p_t$. We want to show that $m\in \Union_{l=1}^jA_{il}$. This then proves that $\Union_{l=1}^jA_{il}$ is poset ideal.

Since $t\in A_{ij}$ and since $A_{ij}\subset \Union_{l=1}^{j-1} A_{i-1,l}$, it follows that there exists $s_{i-1}<s_i=j$ such that $t\in A_{i-1,s_{i-1}}$ Proceeding in this way we find a sequence $s_1<s_2<\cdots< s_i$ such that $t\in \Sect_{l=1}^i A_{l,s_l}$. By the definition of the sorting operator, there must exist a sequence $s_1'<s_2'<\cdots <s_{i-1}'$ with $s_l\leq s_l'\leq s_{l+1}$ for $l=1,\ldots,i-1$ such that $t\in \Sect_{l=1}^{i-1} B_{l,s'_l}$.

Since $t\in A_{ls}$ for some $l$, if and only if $x_{st}|u$, and similarly $t\in B_{ls}$ for some $l$, if and only if $x_{st}|v$  (see (\ref{divides})), it follows that
\[
\prod_{l=1}^ix_{s_lt}\prod_{l=1}^{i-1}x_{s_l't}\quad\text{divides}\quad uv.
\]
Since $uv=u_\MI u_\MJ$ there exists $C$ and $D$ such that
\[
\prod_{c\in C}x_{ct}|u_\MI,\quad \prod_{d\in D}x_{dt}|u_\MJ\quad \text{and}\quad \prod_{c\in C}x_{ct}\prod_{d\in D}x_{dt}=\prod_{l=1}^ix_{s_lt}\prod_{l=1}^{i-1}x_{s_l't}
\]
with $c,d\leq j$ for all $c\in C$ and $d\in D$.

Applying Lemma~\ref{tonight} we conclude that there exist sets $C'$ and $D'$, and bijections $C\to C'$ with $c\mapsto c'\leq c$ and $D\to D'$ with $d\mapsto d'\leq d$ such that
\[
\prod_{c'\in C'}x_{c'm}|u_\MI,\quad \prod_{d'\in D'}x_{d'm}|u_\MJ.
\]
It is clear that $c',d'\leq j$ for $c'\in C'$ and $d'\in D'$.

It follows from the definition of the sorting operator that $i$ of the factors of $\prod_{c'\in C'}x_{c'm}\prod_{d'\in D'}x_{d'm}$ appear in $u$ and $i-1$ of them in $v$. Therefore, by (\ref{divides}) and statement (1) of  Lemma~\ref{terrible}  there exist pairwise different integers  $l_1,\ldots,l_i$ and integers $c_1,\ldots, c_i\leq j$ such that $m\in A_{l_a, c_a}$ for $a=1,\ldots,i$. Therefore there  is at least one $a$ such that $l_a\geq i$. By using part (2) of Lemma~\ref{terrible} we see that $m\in A_{ic}$ for some $c<c_a\leq j$, as desired.

In the same way one shows that $\Union_{l=1}^jB_{il}$ is a poset ideal for all $i$ and $j$.
\end{proof}

\section{The Rees algebra of $H_{r,s}(P)$}

In this section we study the Rees algebra of $H_{r,s}(P)$. Before stating our main result we recall some results of \cite{HHM}, since we shall use them to show that the toric ideal of each power of $H_{r,s}(P)$ has again a quadratic Gr\"obner basis.
First we recall some definitions and results on Rees algebra.

Let
$I=(u_1,\ldots,u_m)$ be a monomial ideal in $K[x_1,\ldots,x_n]$ which is generated in one degree.
Let $R=K[y_{1},\ldots,y_{m}]$ and $L$ be the toric ideal of $K[{u_1},\ldots,{u_m}]$ which is the kernel of the surjective homomorphism
\[
\varphi:R\rightarrow K[u_1,\ldots,u_m]
\]
 defined by $\varphi(y_{i})=u_i$ for all $i$.

\medskip

Let $T$ be the polynomial ring over $K[x_1,\ldots,x_n]$ in the variables
$y_1,\ldots , y_m$. We  may regard  $T$ a bigraded $K$-algebra by setting $\deg(x_i) = (1, 0)$ for
$i = 1,\ldots, n$ and $\deg(y_j) = (0,1)$ for $j = 1,\ldots ,m$.
For any two vectors  $a=(a_1,\ldots,a_n)$ and $b=(b_1,\ldots,b_m)$ with all $0 \leq a_i,b_j$ in $\mathbb{{Z}}$ we write ${\textbf{x}}^a$ for the monomial $x_{1}^{a_1}\cdots x_{n}^{a_n}$ and ${\textbf{y}}^b$ for the monomial $y_{1}^{b_1}\cdots y_{m}^{b_m}$.

Let $\prec$ be a monomial order on $R$. A monomial ${\textbf{y}}^a$ in $R$ is called
a {\em standard monomial} of $L$ with respect to $\prec$, if it does not belong to the initial ideal of $L$. We recall the $\ell$-{\em exchange
property} which was introduced in  \cite{HHM}:

\medskip

The ideal $I$ satisfies the $\ell$-exchange
property with respect to the monomial order $\prec$ on $R$,  if for any two standard monomials $\textbf{y}^a$ and $\textbf{y}^b$ in $L$ of  same degree satisfying
\begin{enumerate}
\item[(i)]
$\deg_{x_t} \varphi(\textbf{y}^a)=\deg_{x_t}\varphi(\textbf{y}^b)$ for $t=1,\ldots,q-1$ with $q\leq n-1$,
\item [(ii)]
$\deg_{x_q}\varphi(\textbf{y}^a)<\deg_{x_q}\varphi(\textbf{y}^b)$,
\end{enumerate}
there exists a factor $u_\delta$ of $\varphi(\textbf{y}^a)$ and $q<j\leq n$ such that  $x_qu_{\delta}/x_j\in I$.

\medskip

The following result is a slight generalization of \cite[Theorem 4.3]{HHM}.

\begin{Proposition}
\label{exchange}
Let $I\subset K[x_1,\ldots,x_n]$ be a weakly polymatroidal ideal which is sortable. Then $I$ satisfies the $\ell$-exchange property with respect to the sorting order.
\end{Proposition}
\begin{proof}

Let $\textbf{y}^a$ and $\textbf{y}^b$ be two standard monomials in $L$ satisfying (i) and (ii). Suppose that $\varphi(\textbf{y}^a)=u_{{i_1}}\cdots u_{{i_d}}$ and $\varphi(\textbf{y}^b)=u_{{j_1}}\cdots u_{{j_d}}$, and that all pairs $(u_{i_l},u_{i_{l'}})$ and $(u_{j_l},u_{j_{l'}})$ are sorted. It follows from (i) that
$\deg_{x_t}(u_{i_l})=\deg_{x_t}(u_{j_l})$ for $l=1,\ldots,d$ and for $t=1,\ldots,q-1$, and (ii) implies that
there exists $1\leq l\leq d$ with $\deg_{x_q}(u_{i_l})<\deg_{x_q}(u_{j_l})$. Since $\deg_{x_t}(u_{i_l})=\deg_{x_t}(u_{j_l})$ for $t=1,\ldots,q-1$ and $\deg_{x_q}(u_{i_l})<\deg_{x_q}(u_{j_l})$, and since $I$ weakly polymatroidal there exists
$j>q$ with $x_qu_{i_l}/x_j\in I$, as desired.
\end{proof}

\medskip

Let $t$ be a variable over
$K[x_1,\ldots,x_n]$. Then the {\em Rees ring}
\[
\MR(I)= \bigoplus_{j=0}^{\infty} I^jt^j\subset K[x_1,\ldots,x_n,t]
\]
is a bigraded algebra with $\deg(x_i) = (1, 0)$ for
$i = 1,\ldots, n$ and $\deg(u_jt) = (0, 1)$ for $j = 1,\ldots,m$. We recall that the toric ideal of $\Rees(I)$ is the ideal $P_{\Rees(I)}\subset T$ which is the kernel of the surjective homomorphism $\varphi: T\rightarrow \Rees(I)$ with  $x_i\mapsto x_i$ for all $i$ and $y_{j}\mapsto u_jt$ for all $j$.

Let $<^{\sharp}$ be an arbitrary monomial order on $R$ and $<_{\lex}$ the lexicographic order on $K[x_1,\ldots,x_n]$ with respect to $x_1>\cdots>x_n$. The new monomial order $<_{\lex}^{\sharp}$ is defined on $T$ as follows: For two monomials  ${\textbf{x}}^a{\textbf{y}}^b$ and  ${\textbf{x}}^{a'}{\textbf{y}}^{b'}$ in $T$, we have
${\textbf{x}}^a{\textbf{y}}^b <_{\lex}^{\sharp}  {\textbf{x}}^{a'}{\textbf{y}}^{b'}$ if and only if
(i) ${\textbf{x}}^a<_{\lex}{\textbf{x}}^{a'}$ or
(ii) ${\textbf{x}}^a={\textbf{x}}^{a'}$ and ${\textbf{y}}^{b}<^{\sharp}{\textbf{y}}^{b'}$.

\medskip
Let $H_{r,s}(P)\subset S$ be generated by the monomials in $G_{r,s}(P)=\{u_1,\cdots,u_m\}$. We have shown in Theorem~\ref{irscm} that $H_{r,s}(P)$ is weakly polymatroidal for the following order of the monomials
\[
x_{11}>x_{21}>\cdots>x_{r1}>x_{12}>\cdots>x_{r2}>\cdots>x_{1n}>\cdots>x_{rn}.
\]
With respect to the same order of the variables the set of monomials is sortable, as shown in Theorem~\ref{fatemehsortable}. Thus if we let $<^{\sharp}$  be the  monomial order given by the property that  $G_{r,s}(P)$ is sortable, we may  apply  \cite[Theorem~5.1]{HHM} to  obtain

\begin{Theorem}
The reduced Gr\"obner basis of the toric ideal $P_{\Rees(H_{r,s}(P))}$ with respect to
$<^{\sharp}$ consists of all binomials belonging to
$\mathcal{G}_{<^{\sharp}}(L_{r,s})$ together with the binomials
\[
x_{ir}y_{k}-x_{js}y_{l},
\]
where $x_{ir}>x_{js}$  with $x_{ir}u_k=x_{js}u_l$ and $x_{js}$ is the largest monomial for which $x_{ir}u_k/x_{js}$ belongs to $H_{r,s}(P)$.
\end{Theorem}

\begin{Corollary}
\label{normalcm}
The Rees ring $\MR(H_{r,s}(P))$ is a normal Cohen--Macaulay domain. In particular all powers of $H_{r,s}(P)$ are normal.
\end{Corollary}

\begin{proof}
We see from the description of $\mathcal{G}_{<_{\lex}^{\sharp}}(P_{\Rees(H_{r,s}(P))})$ that the initial ideal of $P_{\Rees(H_{r,s}(P))}$ is squarefree.  Therefore the result follows from  \cite[Corollary 8.8]{St}  together with  \cite[Theorem 1]{Ho}.
\end{proof}

\begin{Corollary}
\label{Ralf}
The Rees ring $\MR(H_{r,s}(P))$ is Koszul.
\end{Corollary}

\begin{proof}
Since  the initial ideal of  $P_{\Rees(H_{r,s}(P))}$ is generated in  degree 2, the assertion follows from a theorem of Backelin and
Fr\"oberg \cite[Theorem 4(b)]{Ralf}.
\end{proof}

\begin{Corollary}
\label{Ralf}
 All powers of $H_{r,s}(P)$ have a linear resolution.
\end{Corollary}

\begin{proof}
Since for all monomials $u$ in the  minimal set of generators of the initial ideal of  $P_{\Rees(H_{r,s}(P))}$ we have $\deg_{x_{ij}}u\leq 1$ for all $i,j$, the so-called $x$-condition for $H_{r,s}(P)$ is satisfied. Thus  \cite[Corollary 1.2]{HHZ}  yields the desired conclusion.
\end{proof}

As an extension \cite[Corollary 3.8]{HH05} we have

\begin{Corollary}
\label{limit}
Let $P$ be a poset of cardinality $n$. Then
\[
\lim_{k\to \infty}\depth(S/I_r(P)^k)=n-1.
\]
\end{Corollary}

\begin{proof}It is shown in \cite[Proposition 10.3.2]{book} (see also \cite{EH} in combination with \cite[Proposition 1.1]{H}) that
\[\lim_{k\to \infty}\depth(S/I_r(P)^k)=\dim S-\ell(I_r(P)),
 \]
if $\MR(I_r(P))$  is Cohen--Macaulay. Here  $\ell(I_r(P))$ denotes the analytic spread of the ideal $I_r(P)$. Since all generators of $I_r(P)$ have the same degree,  we have that $\ell(I_r(P))=
\dim R_r(P)$. Thus the desired conclusion follows from Corollary~\ref{dim}.
\end{proof}


\begin{thebibliography}{10}
\bibitem{Ralf}
J.~Backelin and R.~Fr{\"o}berg, {\em Koszul algebras, Veronese subrings and rings with linear resolution.} Rev.\
Roumaine Math.\ Pures Appl. 30, 85--97 (1985).

\bibitem{ER}
J.A.~Eagon and V.~Reiner, {\em Resolutions of Stanley-Reisner rings and Alexander duality.} J.\ of Pure and Appl. Algebra, 130, 265 -- 275 (1998).

\bibitem{EH}
D.~Eisenbud and C.~Huneke, {\em Cohen--Macaulay Rees algebras and their specialization} J.\ Algebra, 81,
202--224 (1983).

\bibitem{Ralf1}
R.~Fr{\"o}berg, {\em Koszul algebras,} in ``Advances in Commutative Ring Theory" (D. E. Dobbs, M. Fontana
and S. E. Kabbaj, Eds.), Lecture Notes in Pure and Appl.Math., Volume 205, Dekker, New York, NY,  337--350
(1999).

\bibitem{Villa2} H.~T.~Ha, S.~Morey, E.~Reyes, R.~Villarreal, {\em Cohen-Macaulay admissible clutters.} J. Commut. Algebra 1, no. 1, 463--480   (2009).

\bibitem{cm}
J.~Herzog  and T.~Hibi, {\em Distributive lattices, bipartite graphs
and Alexander duality.}  J. Algebraic Combin.  22,   289--302 (2005).


\bibitem{HH05}
J.~Herzog and T.~Hibi, {\em The depth of powers of an ideal.} J.~Algebra, 291,  534--550 (2005).


\bibitem{book}
J.~Herzog and T.~Hibi, Monomial Ideals. Springer (2010).


\bibitem{HHM}
J.~Herzog, T.~Hibi and M.~Vl\u{a}doiu, {\em Ideals of fibers type and polymatroidals.}
Osaka J.\ Math. 42, 807--829  (2005).


\bibitem{HHZ}
J.~Herzog,  T.~Hibi and X.~Zheng, {\em Monomial ideals whose powers have a
linear resolution.} Math.\ Scand., 95, 23--32 (2004).


\bibitem {HT} J.~Herzog and Y.~Takayama, {\em Resolutions by mapping
cone}. Homology Homotopy Appl. 4, 277--294 (2002).


\bibitem{Hibi1}
T.~Hibi, {\em Distributive lattices, affine semigroup rings and
algebras with straightening laws, in ``Commutative Algebra and
Combinatorics" (M. Nagata and H. Matsumura, eds.)} Adv. Stud. Pure
Math. 11, North-Holland, Amsterdam,  93--109 (1987).


\bibitem{Ho} M.~Hochster,  {\em Rings of invariants of tori, Cohen-Macaulay rings generated by monomials, and polytopes.} Ann.\ of Math.  96, 318--337 (1972).


\bibitem {H}
 C.~Huneke, {\em On the  associated  graded ring of an ideal}. Illinois J.\ Math.
26, 121--137 (1982).


\bibitem{Hibi2}
M.~Kokubo and T.~Hibi, {\em Weakly polymatroidal ideals.} Algebra
Colloq. 13,  711--720  (2006).

\bibitem{Villa1} 
S.~Morey, E.~Reyes, R.H.~ Villarreal, {\em Cohen--Macaulay, shellable and unmixed clutters with a perfect
matching of K\"onig type.} J.\ Pure and Appl.\ Alg. 212,   1770--1786, (2008)

\bibitem{Sta}
R.~Stanley, {\em  Balanced Cohen-Macaulay complexes.} Trans. Amer.\ Math.\ Soc.\ 249, 139--157  (1979).


\bibitem{St}
B.~Sturmfels,  {\em Gr\"obner Bases and Convex Polytopes}. Amer.\ Math.\ Soc., Providence, RI, (1995).


\bibitem{Vi}
R.H.~Villarreal,  Monomial algebras. Marcel Dekker (2001).
\end{thebibliography}
\end{document}